\documentclass[12pt]{amsart}

\usepackage{amsfonts,amssymb,amsmath,graphicx}

\usepackage{color}
\usepackage{amssymb}
\usepackage{hyperref}

\renewcommand{\phi}{\varphi}
\renewcommand{\kappa}{\varkappa}
\renewcommand{\epsilon}{\varepsilon}
\renewcommand{\emptyset}{\varnothing}
\newcommand{\C}{\mathbb C}
\renewcommand{\P}{\mathbb P}
\newcommand{\A}{\mathbb A}
\newcommand{\F}{\mathbb F}
\newcommand{\Z}{\mathbb{Z}}
\newcommand{\R}{\mathbb{R}}
\newcommand{\Q}{\mathbb{Q}}
\newcommand{\bL}{\mathbb L}
\newcommand{\kk}{\mathbf{k}}

\newcommand{\bE}{\mathbf E}

\newcommand{\divisor}{\mathcal D}
\DeclareMathOperator{\FNF}{FNF}
\newcommand{\HH}{\mathbb H}

\DeclareMathOperator{\Mot}{Mot}
\DeclareMathOperator{\sMot}{\overline{Mot}^{\lambda-\mathrm{ring}}}

\DeclareMathOperator{\Spec}{Spec}
\DeclareMathOperator{\cl}{cl}

\DeclareMathOperator{\res}{res}

\DeclareMathOperator{\rk}{rk}

\DeclareMathOperator{\Exp}{Exp}
\DeclareMathOperator{\Log}{Log}

\newcommand{\gl}{\mathfrak{gl}}

\newcommand{\cC}{\mathcal C}

\newcommand{\cF}{\mathcal F}

\newcommand{\cM}{\mathcal M}

\newcommand{\cO}{\mathcal O}
\newcommand{\cP}{\mathcal P}

\newcommand{\cY}{\mathcal Y}

\newcommand{\Bun}{\mathcal{B}un}

\newcommand{\Pair}{\mathcal{P}air}

\newcommand{\Conn}{\mathcal{C}onn}
\newcommand{\Higgs}{\mathcal{H}iggs}

\DeclareMathOperator{\Jac}{Jac}

\newtheorem{theorem}{Theorem}[subsection]
\newtheorem{definition}[theorem]{Definition}
\newtheorem{lemma}[theorem]{Lemma}

\newtheorem{conjecture}[theorem]{Conjecture}

\author{Roman~Fedorov}
\address{Roman Fedorov, University of Pittsburgh, Pittsburgh, PA, USA}
\email{fedorov@pitt.edu}

\author{Alexander Soibelman}
\address{Alexander Soibelman, IHES, 35 route de Chartres, Bures-sur-Yvette, F-91440, France}
\email{asoibel@ihes.fr}

\author{Yan Soibelman}
\address{Yan Soibelman, Kansas State University, Manhattan, KS, USA}
\email{soibel@math.ksu.edu}

\title[Motivic invariants of irregular connections]{Motivic invariants of moduli stacks of Higgs bundles and bundles with  connections:  results and speculations}

\begin{document}
\begin{abstract}
We review some results and techniques from our papers~\cite{FedorovSoibelmans,FedorovSoibelmansParabolic,FedorovSoibelmansIrregular}, devoted to the computation of motivic classes of stacks of parabolic Higgs budles and bundles with connections on a curve.  In the last section we present some directions for future work, as well as some speculations. The latter include a generalization of the  P=W conjecture inspired by the work of Maxim Kontsevich and the third author on the Riemann--Hilbert correspondence for complex symplectic manifolds as well as our running project on the motivic classes of the moduli stacks of nilpotent pairs on the formal disk and geometric Satake correspondence for double affine Grassmannians.
\end{abstract}

\maketitle

\dedicatory{Dedicated to Maxim Kontsevich on occasion of his 60th anniversary.}

\begin{center}
    \includegraphics[width=0.4\textwidth]{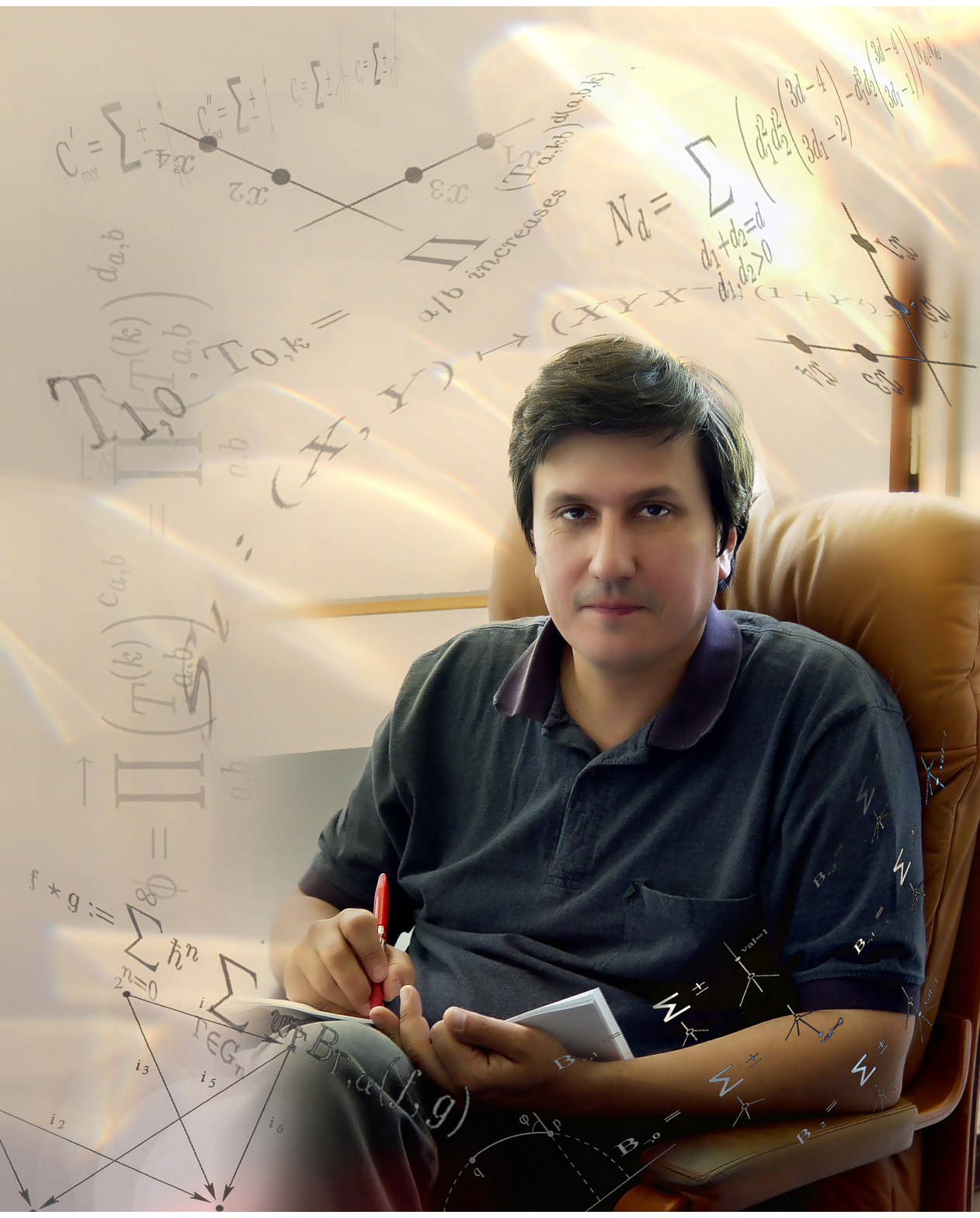}
\end{center}

Some of Maxim's achievements are visible on the photo\footnote{We thank the remarkable photographer and journalist Alla Putintseva for this photo collage.} : Gromov--Witten theory, deformation quantization, Deligne's conjecture, etc. There is also a fragment of the wall-crossing formula from the joint paper with the third author; a categorical version of this formula plays an important part in the current paper.

\tableofcontents

\section{Introduction} The paper is based on the talk given by the third author at IHES in September 2024, at the conference ``Crossroads of mathematics in 20th century'', in honor of Maxim Kontsevich's 60th anniversary. The talk was, in turn, based on our three papers~\cite{FedorovSoibelmans,FedorovSoibelmansParabolic,FedorovSoibelmansIrregular}. In the last paper of the sequence~\cite{FedorovSoibelmansIrregular} we computed motivic Donaldson--Thomas invariants of stacks of irregular semistable parabolic $GL(n)$-Higgs bundles on a smooth geometrically connected projective curve $X$ over a field $\kk$ of characteristic zero, provided that the singular parts of Higgs fields at irregular singular points are generic. We also computed in loc.~cit.\ the motivic Donaldson--Thomas invariants of stacks of parabolic bundles with irregular connections under the same assumptions. We recall that \emph{motivic Donaldson--Thomas invariants\/} were introduced in the joint paper of Maxim Kontsevich and the third author (see~\cite{KontsevichSoibelman08}); an alternative approach was proposed in~\cite{KontsevichSoibelman10}).

\subsection{Acknowledgements}
The work of R.F.~was partially supported by NSF grant DMS-2402553.   The work of A. S. was partially supported by the IHES. The work of Y.S. was partially supported by the  NSF and the Simons Foundation grants as well
as by the IHES Sabbatical Professorship.  A.S. and Y.S. thank  IHES for the excellent research conditions.

\section{Motivic classes of moduli stacks} We are interested in ``motivic volumes'' of the moduli stacks of Higgs bundles and of bundles with connections. We compute the motivic volumes of the moduli stacks of semistable objects of the corresponding categories. More precisely, we fix discrete parameters (in our case, ranks and parabolic degrees) and consider the corresponding moduli stacks of semistable objects. Generating functions of the motivic volumes over the monoid of discrete parameters belonging to a strict cone is an example of a Donaldson--Thomas series introduced by Maxim Kontsevich and the third author in~\cite{KontsevichSoibelman08}.

By definition, these generating functions live in (completed) \emph{motivic Hall algebras}. There is a homomorphism from the latter to the \emph{motivic quantum torus\/} associated with the lattice of Chern classes. Both the motivic Hall algebra and the motivic quantum torus are associative algebras over the commutative ring of \emph{motivic functions}. In this paper, we will only need the ring of motivic functions over a point, which we call the \emph{ring of motivic classes}. We will now recall some basic facts about the latter.

\subsection{Ring of motivic classes} The map $Y\mapsto|Y({\F}_q)|$  taking an algebraic variety $Y$ to its number of points over a finite field is an example of $\Z$-valued \emph{motivic measure}: for any varieties $Y_1$ and $Y_2$ over $\F_q$ we have

(i) $\mu(Y_1\sqcup Y_2)=\mu(Y_1)+\mu(Y_2)$, and

(ii) $\mu(Y_1\times Y_2)=\mu(Y_1)\mu(Y_2)$.

Let $\kk$ be a field. There is a ``universal'' ring $K_0(Var_\kk)$ such that all motivic measures factor through it. As an abelian group it is generated by isomorphism classes $[X]$ of varieties over the ground field $\kk$, subject to the scissor relations
\[
    [X]=[X\backslash Z]+[Z],
\]
for any closed subvariety $Z\subset X$. The commutative ring structure on $K_0(Var_\kk)$ is given by:
\[
 [X]\cdot[Y] = [X \times_\kk Y]_{red}.
\]
The unit element of this ring is $[\Spec\kk]$. The element $[X]\in K_0(Var_\kk)$ is called the \emph{motivic class\/} of $X$. It is easy to see that the motivic class of a vector bundle is the product of the motivic class of the base and the motivic class of the fiber. We have the following examples of motivic measures for $\kk=\C$ (a.k.a ``realizations'' of $K_0(Var_\kk)$):

    $\bullet$ Euler characteristic with compact support: $\chi_c\colon K_0(Var_\C)\to\Z$ sending the class of a variety $X$ to $\sum_{i}(-1)^i \dim H_c^i(X,\Q)$.

    $\bullet$ Deligne--Hodge polynomial $K_0(Var_\C)\to\Z[z_1,z_2]$ sending the class of a variety $X$ to
    \[
        \sum_{i\ge0}(-1)^i\sum_{p,\,q\ge0}\dim Gr^p_F(Gr^W_{p+q} H^i_{DR,c}(X))z_1^{p}z_2^q,
    \]
    where $Gr^W_\bullet$ and $Gr_F^\bullet$ denote the graded components with respect to the weight and Hodge filtration.

Many motivic classes can be computed in terms of the Lefschetz motive $\bL=[\A^1_\kk]$ by decomposing a variety into the union of cells or making it into a vector bundle over a base with known motivic class. For example, we have $[\P^1_\kk]=\bL+1$ and $[GL(n)]=(\bL^n-\bL^{n-1})\ldots(\bL^n-1)$. If $\kk=\F_q$, then applying the realization $[X]\mapsto |X({\F_q})|$  we obtain the well-known formula for the number of points of $GL(n,\F_q)$.

\subsection{Motivic classes of stacks}
One can extend the ring of motivic classes to incorporate motivic classes of Artin (a.k.a.~algebraic) stacks of finite type over $\kk$ with affine stabilizers. Such an algebraic stack $\cY$ can be stratified by global quotients as $\cY=\sqcup_{i=1}^nY_i/GL(n_i)$, where $Y_i$ are varieties. Then one defines
\[
    [\cY]=\sum_{i=1}^n[Y_i]/[GL(n_i)].
\]
This expression does not make sense in $K_0(Var_\kk)$, so we pass to a certain completion $\hat K_0(Var_\kk)$ of $K_0(Var_\kk)[\frac1\bL]$, in which $\frac1\bL$ is treated as a small parameter. In the completed motivic ring $[GL(n)]$ has an inverse.

There are several versions of the notion of motivic ring. We are going to denote by $\Mot(\kk)$ the one that suits our particular needs (note that it is denoted by~$\sMot(\kk)$ in~\cite{FedorovSoibelmansIrregular}). In particular we will need a version that carries a $\lambda$-ring structure. The latter can be  obtained from $\hat K_0(Var_\kk)$ by a certain quotient construction, since $\hat K_0(Var_\kk)$ itself has a pre-$\lambda$-ring structure only (see e.g.~\cite[Sect.~1.8]{FedorovSoibelmansIrregular} where we revisited this well-known constructions in the way most convenient for our purposes). Having the $\lambda$-ring structure allows us to use the plethystic calculus, including the plethystic exponent $\Exp$ and its inverse, the plethystic logarithm $\Log$, as in the usual ring of symmetric functions (see~\cite{macdonald1998symmetric}).

\subsection{Two examples of motivic classes of stacks} The following two results are motivic versions of the calculations of the number of points of the corresponding stacks over finite fields.

\subsubsection{Motivic classes of moduli stacks of bundles on a curve (after Behrend and Dhillon)}
Let $X$ be a smooth projective curve of genus $g$ over a field $\kk$ of characteristic zero.
\begin{theorem}[\cite{BehrendDhillon,FedorovSoibelmans}]
We have in $Mot(\kk)$ the following formula for the motivic class of the stack of vector bundles of rank $r$ and degree $d$:
\[
    [\Bun_{r,d}(X)]=\bL^{(r^2-1)(g-1)}\frac{[\Jac(X)]}{\bL-1}\prod_{i=2}^r\zeta_X(\bL^{-i}),
\]
where $\zeta_X(z)$ is Kapranov's motivic zeta function defined as
\[
    \zeta_X(z):=\sum_{n=0}^{\infty}[Sym^n X]z^n,
\]
while $\Jac(X)$ is the Jacobian of $X$. Here $Sym^nX$ is the $n$-th symmetric power of $X$.
\end{theorem}
This is a motivic analog of Siegel's formula for the number of points over a finite field.

\subsubsection{Motivic G\"ottsche formula}
Let $S/\kk$ be a smooth surface, where  $\kk$ is an algebraically closed field of characteristic zero. Then
\[
    \sum_{n\ge 0}[Hilb_n(S)]z^n=\prod_{m\ge 1}\zeta_S(\bL^{m-1}z^m),
\]
where $Hilb_n(S)$ is the Hilbert scheme of $S$.

\section{For which stacks are we interested in computing motivic classes?}
\subsection{Higgs Bundles and \texorpdfstring{$\epsilon$}{}-connections}
From now on, we fix a field $\kk$ of characteristic zero and a smooth projective, geometrically connected $\kk$-curve $X$. Let $D\subset X(\kk)$ be a finite set of $\kk$-rational points of $X$. Let $\divisor=\sum_{x\in D}n_xx$ be a divisor on $X$, where~$n_x$ are positive integer numbers. To be able to talk about Higgs bundles and bundles with connections on $X$ simultaneously, we introduce the following notion.

\begin{definition}
    Let $\epsilon\in\kk$. \emph{An $\epsilon$-connection\/} on a vector bundle $E$ on $X$ with poles bounded by $\divisor$ is a morphism of sheaves of vector spaces $\nabla:E\to E\otimes\Omega_X(\divisor)$ satisfying the Leibniz rule:
\[
    \nabla(fs)=f\nabla(s)+\epsilon s\otimes df.
\]
Here $f$ and $s$ are any local sections of $\cO_X$ and $E$ respectively.
\end{definition}

If $\epsilon\ne0$ and $\nabla$ is the $\epsilon$-connection, then $\epsilon^{-1}\nabla$ is a singular connection on $E$. If $\epsilon=0$, then the $\epsilon$-connection is the same as a singular Higgs field. We call pairs $(E,\nabla)$ as above ``$\epsilon$-connections'' for brevity. Note also that the curvature of an $\epsilon$-connection is equal to zero as $\dim X=1$. The rank and the degree of $(E,\nabla)$ are by definition the rank and the degree of $E$.

\subsection{Parabolic structures}
For a rational point $x\in X(\kk)$, we denote by $nx$ be the $n$-th infinitesimal neighborhood of $x$. Its ring of functions $k[nx]$ is isomorphic to $\kk[z]/z^n$.
\begin{definition}
Let $E$ be a vector bundle on $X$ and $x\in X(\kk)$. For $n\in\Z_{>0}$, \emph{a level $n$ parabolic structure at $x$ on $E$\/} is a filtration
\begin{equation}
    E|_{nx}=E_0\supset E_1\supset\ldots\supset E_N=0\text{ for }N\gg0
\end{equation}
of the free $k[nx]$-module $E|_{nx}$ by \emph{free\/} submodules. For a divisor $\divisor=\sum_{x\in D}n_xx$ with $n_x>0$, \emph{a level~$\divisor$ parabolic structure on $E$\/} is a collection $E_{\bullet,\bullet}=(E_{x,j})$ where $E_{x,j}$ is a level $n_x$ parabolic structure on~$E$ at $x$ for all $x\in D$. The pair $(E,E_{\bullet,\bullet})$ is called a \emph{level $\divisor$ parabolic bundle on $X$}.
\end{definition}

Given a level $n$ parabolic structure at $x$ on $E$, one can always trivialize $E$ near $x$ so that each $E_i$ is the submodule of $E|_{nx}$ generated by $e_1,\ldots,e_{\rk E_i}$, where $e_1,\ldots,e_{\rk E}$ is the basis of $E|_{nx}$ provided by the trivialization. Such a trivialization is called \emph{compatible with the parabolic structure}.

\subsection{Stack of parabolic irregular \texorpdfstring{$\epsilon$}{}-connections with prescribed normal forms}
Let us describe the stacks we are dealing with in a more precise way.

Recall that $\divisor$ is an effective divisor on $X$ with set-theoretic support $D\subset X(\kk)$. We define the set of  \emph{full normal forms} of future $\epsilon$-connections as
 \[
    \FNF(\divisor):=\prod_{x\in D}(\Omega_X(n_xx)/\Omega_X)^{\Z_{>0}}.
\]
Let $\epsilon\in\kk$ and $\zeta=(\zeta_{x,j}\,|\;x\in D,j\in\Z_{>0})\in\FNF(\divisor)$. Then $\Conn(\epsilon,X,\divisor,\zeta)$ is the {moduli stack} classifying triples $(E,E_{\bullet,\bullet},\nabla)$, where\\
    $\bullet$ $(E,\nabla)$ is an $\epsilon$-connection on $X$ with poles bounded by $\divisor$;\\
    $\bullet$ $E_{\bullet,\bullet}=(E_{x,j})$ is a level $\divisor$ parabolic structure on $E$;\\
    $\bullet$ for $x\in D$ and a local trivialization of $E$ near $x$ compatible with the parabolic structure write
    \[
    \nabla=\epsilon d+A
    \]
    and view $A$ as a block matrix with blocks of sizes $\rk(E_{x,i-1}/E_{x,i})\times\rk(E_{x,j-1}/E_{x,j})$. Then the polar part of the $i$-th diagonal block is equal to $\zeta_{x,i}$, and the blocks lying below the diagonal have no pole at $x$.

Consider the monoid
\[
    \Gamma_D=\{\gamma=(r,r_{\bullet,\bullet})| r\in\Z_{\ge0}, r_{\bullet,\bullet}\in D\times\Z_{>0}, \forall x\in D, \sum\limits_{j=1}^{\infty}r_{x,j}=r\}.
\]
If $\bE:=(E,E_{\bullet,\bullet})$ is a level $\divisor$ parabolic bundle on $X$, we define its \emph{class\/} by
\[
    \cl(\bE):=(\rk E,(\rk E_{x,j-1}-\rk E_{x,j})_{x\in D,j>0},\deg E)\in\Gamma_D\times\Z.
\]
We have a decomposition according to the classes of parabolic bundles
\[
 \Conn(\epsilon,X,\divisor,\zeta)=\bigsqcup_{\gamma\in\Gamma_D,d\in\Z}\Conn_{\gamma,d}(\epsilon,X,\divisor,\zeta).
\]
We note that the stacks $\Conn_{\gamma,d}(\epsilon,X,\divisor,\zeta)$ are of finite type when $\epsilon\ne0$, see \cite[Cor.~7.2.4]{FedorovSoibelmansIrregular}. One motivation to study these stacks is that they are twisted cotangent bundles of the stack of parabolic vector bundles, see~\cite[Sect.~1.16]{FedorovSoibelmansIrregular}. The corresponding $DT$-series (see~\cite{KontsevichSoibelman08} for the definition) is
\[
 \sum_{\gamma\in\Gamma_D,d\in\Z}[\Conn_{\gamma,d}(\epsilon,X,\divisor,\zeta)]w^{\gamma}z^d\in
 Mot(\kk)[[\Gamma_D,\Z]].
\]

\section{Motivic Hall algebras and quantum tori}
Our results and calculation fall in the broad categorical framework from the paper of Kontsevich and the third author~\cite{KontsevichSoibelman08}, which we briefly recall below.

Let $\cC$ is a $\kk$-linear category such that the stack of objects $\cM=Ob(\cC)=\sqcup_{\gamma\in \Gamma}{\cM}_\gamma$ is a graded ind-Artin (or ind-constructible) stack, where $\Gamma$ is the free abelian group of ``Chern classes" of objects. E.g., $\Gamma=\Z^2$ for the category of coherent sheaves on our curve $X$, and $\gamma=(r,d)$, where $r$ is the rank and $d$ is the degree of a vector bundle. The \emph{motivic Hall algebra $H_\cC$\/} over $Mot(\kk)$ is an associative algebra generated by isomorphism classes of constructible families of objects over $\cM$. In particular it includes delta-functions $\delta_{[E]}:=[E]$ supported at the isomorphism class of an object $E$. The associative product comes from the stack of exact triangles $E\to F\to G$ via the usual pullback-pushforward construction using projections to each of the terms of the exact triangle.

If $\cC$ carries a Bridgeland stability structure with the central charge $Z\colon\Gamma\to\C$, then for any strict sector $V\subset\R^2$ with the vertex at the origin we define a full subcategory $\cC(V)$ generated by semistables $E$ with $Z(E)\in V$. The motivic volume of the stack $Ob(\cC(V))$ is the following expression
\[
    A_V^{Hall}=\sum_{E\in\cC(V)}{\frac{[E]}{[Aut(E)]}}\in\widehat{H}_{\cC(V)}.
\]
The hat means that as an infinite series it belongs to a certain completion of $H_{\cC(V)}$. The Harder--Narasimhan property implies the following factorization (a.k.a wall-crossing) formula
\[
    A_{V}^{Hall}=\prod_{l\subset V}A_l^{Hall},
\]
where the clockwise product in the RHS is taken over all rays $l\subset V$ with the vertex at the origin.

The motivic quantum torus $R_\cC(V)$ is the algebra over $Mot(\kk)$ with generators $e_\gamma,\gamma\in \Gamma$ such that $e_\gamma e_\mu=\bL^{{\frac12}\langle\gamma,\mu\rangle}e_{\gamma+\mu}$, where $\langle\bullet,\bullet\rangle$ is the skew-symmetrization of the Euler form on $K_0(\cC)$. The following is the key result for computations (see~\cite{KontsevichSoibelman08}).

\begin{theorem}
If $\cC$ is an ind-constructible locally ind-Artin 3-di\-men\-si\-onal Calabi--Yau category, then there is a homomorphism of algebras $H_{\cC(V)}\to R_{\cC(V)}$.
\end{theorem}

Applying the homomorphism $H_\cC(V)\to R_\cC(V)$ to the factorization formula in the motivic Hall algebra we obtain the factorization formula in the completion of the motivic quantum torus
\[
    A_V^{mot}=\prod_{l\subset V}A_l^{mot}.
\]
The LHS is called \emph{the motivic DT series} of $\cC(V)$. Take $V=l$. Typically, $H_{\cC(l)}$ is a symmetric algebra of a free $Mot(\kk)$-module. In this case,
\[
    A_l^{mot}=\Exp\left({\sum_{Z(\gamma)\in l}\Omega(\gamma)\frac{e_\gamma}{1-\bL^{-1}}}\right),
\]
where $\Exp$ is the plethystic exponent and $\Omega(\gamma)\in Mot(\kk)$ is \emph{the motivic DT-invariant} of class $\gamma$. Morally, $\Omega(\gamma)$ is the motivic class of the moduli stack of semistables of slope $l$ and Chern class $\gamma$.

Let us make a short historical remark. The homomorphism from the motivic Hall algebra to the quantum torus is sometimes called ``the integration map''. It has a predecessor in the case of finite fields, which was introduced by Reineke in his work on quantum groups. Also, a version of the above homomorphism at the level of Lie algebras (not associative algebras) was introduced by Joyce in his work on the invariants of stacks of semistable objects of \emph{abelian}  categories. All the above papers were not known to the authors of~\cite{KontsevichSoibelman08} at the time of writing. Instead, they developed in~\cite{KontsevichSoibelman08} an approach based on the ideas of motivic integration (in particular, utilizing motivic Milnor fibers) in the framework of what they called ind-constructible Calabi--Yau categories These categories are {\it triangulated}. Later, they proposed an alternative approach in~\cite{KontsevichSoibelman10} by introducing Cohomological Hall algebras. The latter depend on a choice of  $t$-structure.

The starting point for developing the theory of motivic (and numerical) Donaldson--Thomas theory in~\cite{KontsevichSoibelman08} in the categorical framework was the work of the authors on a new type of wall-crossing formulas (see~\cite{KontsevichSoibelman08} for the details), partially motivated by their previous work on Homological Mirror Symmetry.

In order to finish this historical digression, here is a picture, which might be of some historical interest.
It is a photo of Maxim Kontsevich and the third author explaining to Don Zagier the above-mentioned formalism of motivic DT-invariants on a boat trip during the Arbeitstagung--2007 in Bonn.
\begin{center}
\includegraphics[width=1.0\textwidth]{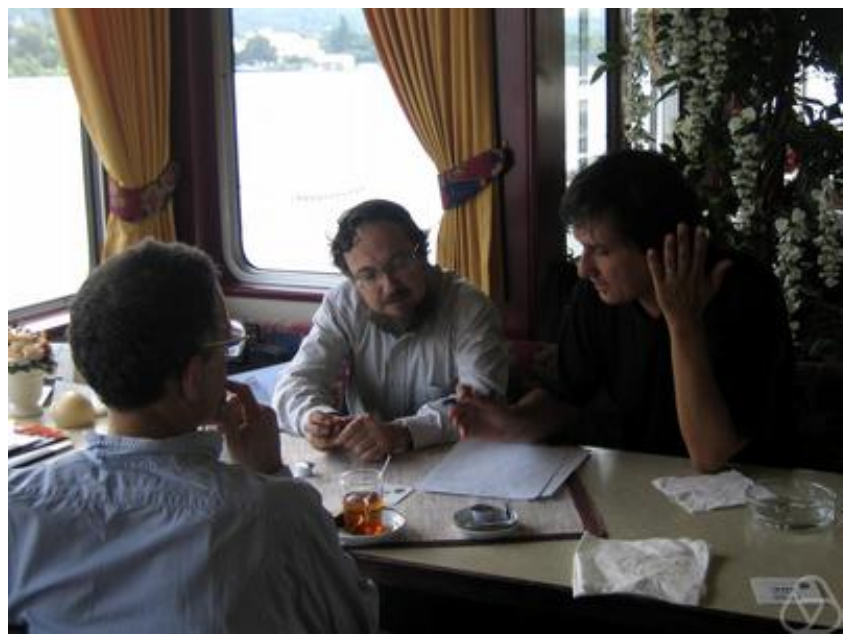}
\end{center}

\subsection{A version for 2-dimensional Calabi--Yau categories}
There is a version of the above formalism for $2d$ Calabi--Yau categories, like e.g. the category $\Higgs_X$ of Higgs bundles on $X$, where $D=\emptyset$. Then, the skew-symmetrized Euler form is trivial, and the motivic quantum torus is commutative, hence isomorphic to $\Mot(\kk)[\Gamma]$. In this case, the motivic $DT$-series becomes a series in commuting variables $w,z$ of the type $\sum_{\gamma, d}a_{\gamma,d}w^{\gamma}z^d$, where $a_{\gamma,d}\in Mot(\kk)$, $(\gamma,d)\in\Gamma$.

In the case of singular Higgs bundles, the lattice $\Gamma=\Z^2$ should be replaced with $\overline\Gamma_D\times\Z$, where $\overline\Gamma_D$ is the free abelian group obtained from the monoid $\Gamma_D$ via the Grothendieck construction. The categories $\cC(V)$ that appear in our story are relatives of the categories of Higgs bundles on $X$. More precisely, it is required that the underlying vector bundles have Harder--Narasimhan spectrum in a strict sector $V$ (e.g., we can take $V=\{(x,y)\in\R^2|x<0,y>0\}$).

\section{Summary of results}
\subsection{Motivic \texorpdfstring{$L$}{}-functions and subsequent computation}

Define the motivic L-function of $X$
\[
 L_X(z):=\zeta_X(z)(1-z)(1-\bL z)\in Mot(\kk)[z],
\]
which is actually a polynomial of degree $2g$. Let $\delta$ be a nonnegative integer. The following series is a generalization of those of Hausel and Mellit
\begin{multline*}
\Omega^{mot}_{X,D,\delta}
=\sum_{\mu\in\cP}w^{|\mu|}(-\bL^\frac12)^{(2g+\delta)\langle\mu,\mu\rangle}\times
\\z^{2\delta n(\mu')}\prod_{\Box\in\mu}\frac{L_X(z^{2a(\Box)+1}\bL^{-l(\Box)-1})}
 {(z^{2a(\Box)+2}-\bL^{l(\Box)})(z^{2a(\Box)}-\bL^{l(\Box)+1})}\times
\\\prod_{x\in D}\tilde H^{mot}_\mu(w_{x,\bullet};z^2).
\end{multline*}
The sum is taken over all partitions and the product is taken over all boxes in the partition. Here $a(\Box)$ and $l(\Box)$ are respectively the arm and the leg lengths of the box, $\mu^\prime=\mu'_1\ge\mu'_2\ge\ldots$ is the conjugate partition, and $n(\mu')=\sum_i(i-1)\mu'_i$. The multiples $\tilde H^{mot}_\mu$ are motivic versions of the modified Macdonald polynomials.

We expand the denominators in positive powers of $z$, so that we have $\Omega^{mot}_{X,D,\delta}\in Mot(\kk)[\bL^{\frac12}][[\Gamma_D,z]]$. In fact, $\Omega^{mot}_{X,D,\delta}$ is symmetric in each of the sequences of variables $w_{x,\bullet}$. Note that the constant term of $\Omega^{mot}_{X,D,\delta}$ is 1, so we can consider
\begin{equation}
 \HH^{mot}_{X,D,\delta}:=(1-z^2)\Log\Omega^{mot}_{X,D,\delta},
\end{equation}
where $\Log$ is the plethystic logarithm.

\begin{theorem}[\cite{MellitNoPunctures,FedorovSoibelmansIrregular}]
Coefficients of $\HH^{mot}_{X,D,\delta}$ are polynomials in $z^{\pm1}$.
\end{theorem}

Our motivic classes will be calculated in terms of $\HH^{mot}_{X,D,\delta}$ but we need some preliminary definitions. Let $\gamma=(r,r_{\bullet,\bullet})\in\Gamma_D$. We will say that $\gamma$ is \emph{full at $x$\/} if $r_{x,j}\in\{0,1\}, \forall j>0$. Let $\cl(E,E_{\bullet,\bullet})=(\gamma,d)$. Then, $\gamma$ is full at $x\in D$ if and only if $E_{x,\bullet}$ is a full flag of submodules of $E_{n_xx}$ (with repetitions, that is, for some $j$ we have $E_{x,j-1}=E_{x,j}$).

Let $\zeta=(\zeta_{x,j})\in\FNF(\divisor)$. We will say that $\zeta$ is \emph{non-resonant for $\gamma$ at $x\in D$\/} if for all $1\le i<j$ such that $r_{x,i}\ne0\ne r_{x,j}$, the polar part $\zeta_{x,i}-\zeta_{x,j}$ has a pole of order exactly $n_x$.

In the following formula, we take $\delta=\delta_\divisor:=\sum_{x\in D}(n_x-1)$. In other words, it is the irregularity index
$\delta=\deg\divisor-|D|\in\Z_{\ge0}$ of $\divisor$, where $|D|$ is the number of points in $D$. For $\gamma=(r,r_{\bullet,\bullet})\in\Gamma_D$ and $\zeta=(\zeta_{x,j})\in\FNF(\divisor)$ set
\[
\gamma\star\zeta:=\sum_{x\in D}\sum_{j>0}r_{x,j}\res\zeta_{x,j}\in\kk.
\]
We also define
\[
    \chi(\gamma)=\chi(\gamma,g,\divisor):=(2g-2)r^2-\delta r+2\sum_{x\in D,i<j}n_xr_{x,i}r_{x,j}.
\]
The following is the first main result of~\cite{FedorovSoibelmansIrregular}.
\begin{theorem}
Assume that $\epsilon\ne0$, the class $\gamma'\in\Gamma_D$ is full at all $x\in D$ such that $n_x\ge2$, and $\zeta\in\FNF(\divisor)$ is non-resonant for $\gamma'$ at all such points $x$. Assume $d'\in\Z$ is such that $\epsilon d'+\gamma'\star\zeta=0$. Then, the motivic class of $\Conn_{\gamma',d'}(\epsilon,X,\divisor,\zeta)$ in $Mot(\kk)$ is equal to the coefficient at $w^{\gamma'}$ in
\[
 (-\bL^\frac12)^{\chi(\gamma')} \Exp\left(\bL\left(\HH^{mot}_{X,D,\delta}\Bigl|_{z=1}\right)_{\gamma\star\zeta\in\epsilon\Z}\right),
\]
where the subscript $\gamma\star\zeta\in\epsilon\Z$ stands for the sum of the monomials whose exponents satisfy the condition.
\end{theorem}
The reader may notice that  by the previous theorem we can evaluate $\HH^{mot}_{X,D,\delta}$ at $z=1$, hence the formula makes sense. We also note that the stack under consideration is empty if $\epsilon d'+\gamma'\star\zeta=0$. Finally, we note that we are only able to calculate the motivic class in the case when the flags are full and the $\epsilon$-connection is non-resonant at all points with $n_x\ge2$. In fact, the parabolic structure at such points is uniquely determined by the $\epsilon$-connection. On the other hand, there are no conditions at the points where $n_x=1$; the parabolic structure is, in general, extra data at such points.

\subsection{Explicit formulas for motivic classes of irregular semistable parabolic \texorpdfstring{$\epsilon$}{}-connections}\label{sect:ExplAnswers2}
In the previous theorem, we calculated the motivic classes of $\epsilon$-connections under the assumption that $\epsilon\ne0$, excluding thus the case of Higgs bundles. In fact, when $\epsilon=0$, the stacks are of infinite type and even of infinite motivic volume. To obtain a finite answer, we introduce the slope stability condition and compute the motivic classes of the semistable loci of $\Conn_{\gamma,d}(\epsilon,X,\divisor,\zeta)$ under the same genericity assumptions.

If $E$ is a vector bundle on $X$, its subbundle $F$ is called  \emph{strict}, if $E/F$ is torsion free. Let $(E,E_{\bullet,\bullet},\nabla)$ be a point of $\Conn_{\gamma,d}(\epsilon,X,\divisor,\zeta)$ and let $F\subset E$ be a strict subbundle preserved by $\nabla$. A new difficulty in the irregular case is that we do not necessarily get an induced parabolic structure on~$F$.

\begin{lemma}[\cite{FedorovSoibelmansIrregular}]
 Let $(E,E_{\bullet,\bullet},\nabla)\in\Conn_{\gamma,d}(\epsilon,X,\divisor,\zeta)$ and let $F\subset E$ be a strict subbundle preserved by $\nabla$. Suppose that $\gamma$ is full at all $x\in D$ such that $n_x\ge2$, and $\zeta$ is non-resonant for $\gamma$ at all such $x$. Then we get an induced parabolic $\epsilon$-connection $(F,F\cap E_{\bullet,\bullet},\nabla|_F)\in\Conn_{\gamma',d'}(\epsilon,X,\divisor,\zeta)$ for some $(\gamma',d')\in\Gamma_D\times\Z$.
\end{lemma}

A sequence $\sigma=(\sigma_{\bullet,\bullet})$ of real numbers indexed by $D\times\Z_{>0}$ is called a \emph{sequence of parabolic weights of type $(X,\divisor)$\/} if for all $x\in D$ we have $\sigma_{x,1}\le\sigma_{x,2}\le\ldots$ and for all $j$ we have $\sigma_{x,j}\le\sigma_{x,1}+n_x$. Set $\gamma\star\sigma:=\sum_{x\in D}\sum_{j>0}r_{x,j}\sigma_{x,j}\in\R$. We call the parabolic $\epsilon$-connection $(E,E_{\bullet,\bullet},\nabla)$ as above \emph{$\sigma$-semistable\/} if for all strict subbundles $F\subset E$ preserved by $\nabla$, we have
\[
 \frac{d'+\gamma'\star\sigma}{\rk F}\le\frac{d+\gamma\star\sigma}{\rk E},
\]
where $(\gamma',d')$ is the class of $(F,F\cap E_{\bullet,\bullet})$. We denote the open substack of $\Conn_{\gamma,d}(\epsilon,X,\divisor,\zeta)$ classifying $\sigma$-semistable parabolic $\epsilon$-connections by $\Conn^{\sigma-ss}_{\gamma,d}(\epsilon,X,\divisor,\zeta)$. One can show that these stacks are always of finite type.

\begin{theorem}[\cite{FedorovSoibelmansIrregular}]
Assume that $\gamma'\in\Gamma_D$ is full at all $x\in D$ such that $n_x\ge2$, and that $\zeta\in\FNF(\divisor)$ is non-resonant for $\gamma'$ at all such $x$. Assume that $d'\in\Z$ is such that $\epsilon d'+\gamma'\star\zeta=0$. Let $\sigma$ be a sequence of parabolic weights of type $(X,\divisor)$. Then, the motivic class of $\Conn^{\sigma-ss}_{\gamma',d'}(\epsilon,X,\divisor,\zeta)$ in $Mot(\kk)$ is equal to the coefficient at $w^{\gamma'}$ in
\[
 (-\bL^\frac12)^{\chi(\gamma')}\Exp\left(\bL\left(\HH^{mot}_{X,D,\delta}\Bigl|_{z=1}\right)_{\substack{\gamma\star\zeta\in\epsilon\Z\\ \gamma\star\sigma-\tau\rk\gamma\in\Z }}\right),
\]
where $\tau:=(d'+\gamma'\star\sigma)/\rk\gamma'$, and $rk(r,r_{\bullet,\bullet}):=r$.
\end{theorem}

\section{Summary of the calculation}
\subsection{Summary of the calculation of the motivic classes in the non-singular case}
In~\cite{FedorovSoibelmans} we followed the approach of~\cite{SchiffmannIndecomposable} and~\cite{MozgovoySchiffmann2020} with modifications necessary for the motivic case. One of the ideas of Schiffmann is to trade semistability for a more manageable condition that still ensures that our stacks are of finite type: we look at Higgs bundles such that the underlying vector bundle is ``nonnegative''. That is, it has no surjective morphisms to a bundle of negative degree. This is equivalent to the Harder--Narasimhan spectrum of the bundle being nonnegative (see~\cite[Sect.~3.2]{FedorovSoibelmans}). We emphasize that this is the condition on the vector bundle rather than on the Higgs pair. We use a combination of a version of the Kontsevich--Soibelman factorization formula and a limiting process to express the motivic classes of semistable Higgs bundles in terms of the motivic classes above (see~\cite[Sect.~3.3 and~3.6]{FedorovSoibelmans}).

To calculate the motivic classes of bundles with connections, we note that for a given vector bundle the space of connections is an affine space modelled on the space of Higgs bundles, provided that the bundle has a connection. We use the theorem of Atiyah claiming that a bundle carries a connection if and only if all its direct summands have degree (equivalently, slope) zero. We then use a product formula (see~\cite[Sect.~3.5]{FedorovSoibelmans}) to calculate the motivic classes of Higgs bundles whose underlying vector bundle is ``isoslopy'', that is, all its direct summands have the same slope. This product formula, while looking very similarly to the one from Kontsevich--Soibelman~\cite{KontsevichSoibelman08}, is closely related to Hua's formula for the number of indecomposable representations of a quiver.

We show that the motivic class of the stack of rank $r$ bundles with connections is equal to the motivic class of the stack of the \emph{semistable\/} Higgs bundles of rank $r$ and degree 0 (see~\cite[Thm.~1.2.1]{FedorovSoibelmans}). This can be called \emph{motivic non-abelian Hodge theory}. The above observation was one of  the motivations for our work. Note that in the twisted case, specificallly, when we look at bundles with connections whose degree and rank are coprime, and the connections are allowed a first order pole with an appropriate scalar residue, a similar statement follows from geometric considerations. In fact, it is shown in \cite[Sect.~4]{HoskinsLehalleurOnVoevodskyMotive} that the Voevodsky motives are equal. However, in the general case we are not aware of a geometric argument. The reason is that the relation between the stack and the coarse moduli space is complicated because there are strictly semistable Higgs bundles. We use this opportunity to emphasize that we are always working with stacks.

We are left with calculating the motivic class of the stack of Higgs bundles with nonnegative underlying vector bundles. Using the Riemann--Roch theorem and Serre's duality, we replace bundles with Higgs fields by bundles with endomorphisms (see~\cite[Lemma~3.5.2]{FedorovSoibelmans}). Using the formalism of power structures and a version of Jordan decomposition, we reduce our calculation to the calculation of the motivic class of the stack of nonnegative vector bundles with \emph{nilpotent\/} endomorphisms (see~\cite[Sect.~3.8]{FedorovSoibelmans}).

Following~\cite{GarciaPradaHeinlothSchmitt} and~\cite{SchiffmannIndecomposable}, we reduce calculating the motivic class of the stack of bundles with nilpotent endomorphisms to calculating motivic classes of chains (that is, filtered vector bundles). We calculate the latter motivic class, following the calculation over finite fields from~\cite{SchiffmannIndecomposable} but using the motivic Hall algebra instead of the usual one (see~\cite[Sect.~5 and~6]{FedorovSoibelmans}). An important part of our calculation is what we call ``motivic Harder's formula''. Intuitively, it says that the motivic class of the stack of full flags of fixed degrees on a vector bundle $E$ is ``approximately'' independent of $E$ when the degrees of quotients tend to $-\infty$  (see~\cite[Sect.~4]{FedorovSoibelmans}). We note also, that these results have been recently extended to finite characteristic, see~\cite{Li2025MotivicClasses}.

\subsection{Summary of the calculation of the motivic classes in the case of regular singularities}
In~\cite{FedorovSoibelmansParabolic} we extended the results of~\cite{FedorovSoibelmans} to the case of parabolic Higgs bundles and connections with regular singularities. Our moduli stacks are stacks of twisted cotangent bundles of the stack of parabolic bundles. The twisting is given by the eigenvalues of the residues and a parameter $\epsilon$ distinguishing Higgs bundles from connections. We also have  stability conditions coming from parabolic weights. The stability conditions and the twisting parameters appear in a symmetric way in the final answer, leading to various ``motivic non-abelian Hodge isomorphisms'', see~\cite[Sect.~9.2]{FedorovSoibelmansParabolic}. We also use the motivic formulas to derive the non-emptiness of various stacks, which is a version of the Deligne--Simpson problem~\cite[Sect.~9.3]{FedorovSoibelmansParabolic}.

The strategy is similar to the non-singular case: we reduce the calculation of our motivic classes to calculating motivic classes of \emph{parabolic\/} bundles with nilpotent endomorphisms, whose underlying vector bundle is nonpositive (for technical reasons we switch from nonnegative to nonpositive bundles). This calculation follows the strategy of~\cite{MellitPunctures} and consists of two parts: first one shows that the motivic class under consideration is the product of a motivic class depending only on the curve (but not on the number of punctures) and the classes attached to the punctures, which are independent of the curve. The former class was calculated in our work~\cite{FedorovSoibelmans} as explained above, while the ``local classes'' are calculated by considering $\P^1$ with two punctures. This calculation can be derived from Mellit's result, or, alternatively from~\cite{Singh2021Counting}, where a similar counting problem was solved for any reductive group $G$ over finite fields (but the calculation is essentially motivic). As for the factorization formula, it is more delicate in the motivic case and requires working carefully with finite and infinite order jets.

\subsection{Summary of our approach in the irregular case}
We now assume that the divisor $\divisor$ is \emph{not\/} reduced. We start by generalizing our question: let $\divisor'$ be an effective divisor such that $\divisor'\le\divisor$. We consider stacks $\Conn^{prtl,-}(\epsilon,X,\divisor,\divisor',\zeta)$ classifying $\epsilon$-connections with \emph{level $\divisor'$\/} parabolic structure and \emph{partially fixed\/} formal normal forms (as before, we work with nonnegative vector bundles).

Our argument is by induction on $\deg\divisor'$: we start with the case when $\divisor'=D$ is a reduced divisor; we recover our previous stacks when $\divisor=\divisor'$.

To run this argument, we had to develop an existence criterion for an $\epsilon$-connection on a parabolic bundle with fully or partially fixed formal normal form, see~\cite[Thm.~1.15.1]{FedorovSoibelmansIrregular}. This criterion implies that for every $\zeta$, every level $D$ parabolic bundle carries an $\epsilon$-connection in $\Conn^{prtl,-}(\epsilon,X,\divisor,D,\zeta)$. This allows to relate the motivic class of the above stack to that of the stack of twisted nilpotent pairs. Twisted nilpotent pairs in our case are pairs $(\bE,\Psi)$, where $\bE=(E,E_{\bullet,\bullet})$ is a level $D$ parabolic bundle and $\Psi\colon E\to E\otimes \cO_X(D-\divisor)$ is a morphism  compatible with the parabolic structure. Nilpotency of the morphism means that iterating it finitely many times we obtain zero. In our case the nilpotency property follows from the observation that the line bundle $\cO_X(D-\divisor)$ is negative. The calculation of the motivic class of the  stack of twisted nilpotent pairs is similar to the regular case considered in~\cite{FedorovSoibelmansParabolic}.

Once the motivic classes of $\Conn^{prtl,-}(\epsilon,X,\divisor,\divisor',\zeta)$ are available for reduced $\divisor'$, we proceed by induction on $\deg\divisor'$ to calculate these motivic classes for any $\divisor'<\divisor$. The idea is to relate the motivic class above to the motivic class of $\Conn^{prtl,-}(\epsilon,X,\divisor,\divisor'-x,\zeta)$, where $x\in D$ is any point with $n'_x\ge2$.

 The case $\divisor'=\divisor$ is more difficult. The idea is to relate the stack $\Conn^-(\epsilon,X,\divisor,\zeta)$ to $\Conn^{prtl,-}(\epsilon,X,\divisor,\divisor-x,\zeta)$. To this end, we show that a parabolic bundle $\bE$ admits an $\epsilon$-connection with formal normal form $\zeta$, if and only if it admits such an $\epsilon$-connection with normal form fixed up to level $\divisor-x$ \emph{and\/} every direct summand of $\bE$ has degree zero.

The rest of the calculation is similar to the regular case. Once the formulas for the motivic classes are obtained, we
apply the technique of Mellit~\cite{MellitIntegrality,MellitNoPunctures,MellitPunctures} to simplify our formulas. Note that this simplification is only possible in the $\lambda$-ring, while the pre-$\lambda$-ring structure is not sufficient. We also note that our results in~\cite{FedorovSoibelmansIrregular} give simpler formulas in the regular and even in the non-singular cases.

Our calculation of motivic classes of stacks of parabolic $\epsilon$-connections is then used to calculate the E-polynomials and the virtual Poincar\'e polynomials of these stacks. As a result, we provide evidence for a conjecture of Diaconescu, Donagi, and Pantev (see~\cite[(1.8)]{DiaconescuDonagiPantev}), which, in turn, is based on a conjecture of Hausel, Mereb, and Wang (see~\cite[Conj.~0.1.1]{HauselMerebWong}).

\section{Final remarks, questions and speculations}
We mention below a few directions for the future work. In the last two subsections we discuss a speculative generalization of the $P=W$ conjecture related to the work of Maxim Kontsevich and the third author on the generalized Riemann--Hilbert correspondence as well as  our running project on motivic classes of nilpotent pairs on the formal disk.

\subsection{Generalization to arbitrary groups}
It would be very interesting to generalize our results and techniques to the moduli stacks of Higgs bundles and bundles with connections for an arbitrary $\kk$-split reductive group $G$. The corresponding categories are non-additive, so techniques related to motivic Hall algebras are not available. One can try to replace motivic Hall algebras by their versions based on the parabolic induction for $G$.

\subsection{Approach via Cohomological Hall algebras}
There is a natural question of how to compute motivic DT-invariants using Cohomological Hall algebras defined in the paper of Maxim Kontsevich and the third author in~\cite{KontsevichSoibelman10}. According to the Section 7.10 of~\cite{KontsevichSoibelman10} the result should agree with our computations obtained via motivic Hall algebras. There is some work on the Cohomological Hall algebras of the moduli stacks of Higgs bundles and Higgs sheaves (see, e.g.,~\cite{SalaSchiffmannCOHAHiggs},~\cite{DPSSV2025COHAYangians}), but they are mostly devoted to  other questions. In any case, we are not aware of the paper in which our results were obtained via Cohomological Hall algebra techniques.

\subsection{Motivic Deligne--Simpson correspondence}
The Deligne-\\ Simpson problem refers to the question of the existence of a Higgs bundle (additive problem) or a bundle with connection (multiplicative problem) on a curve with prescribed singular behavior at given points.

We used motivic techniques in order to formulate the motivic version of the Deligne--Simpson problem (cf.~\cite[Sect.~9.3]{FedorovSoibelmansParabolic}). On the other hand, Jakob--Yun~\cite{JakobYun2023DeligneSimpsonIrregular} proposed an answer in the irregular case in terms of representations of rational Cherednik algebras. It would be interesting to compare their proposal with our computations.

\subsection{Generalized non-abelian Hodge theory and generalized \texorpdfstring{$P=W$}\, conjecture}
Note that Higgs bundles can be interpreted as coherent sheaves on $T^\ast X$ that have pure 1-dimensional support. One can  generalize the problem and ask about motivic classes of moduli stacks of coherent sheaves on  {\it symplectic surfaces} that have pure $1$-dimensional support. This question can be put in the framework of the ongoing project of Maxim Kontsevich and the third author on the \emph{generalized Riemann--Hilbert correspondence\/} and \emph{generalized non-abelian Hodge theory}, in which cotangent bundles  are replaced by arbitrary complex symplectic manifolds. The key idea of the generalized Riemann--Hilbert correspondence is to replace the category of constructible sheaves in the conventional Riemann--Hilbert correspondence by an appropriate Fukaya category. In the case of cotangent bundles both categories agree.

Let us say a few words about the generalized non-abelian Hodge theory in dimension one (details will appear in~\cite{KontsevichSoibelmanTwistor}, see also ~\cite{KontsevichSoibelman2024ExpIntegrals}). In order to save space we will outline below the main ideas without trying to be more precise.

Suppose we are given a complex symplectic surface $(S,\omega^{2,0})$ together with some other choices. The main one is a choice of a partial log-compactification $S_{log}\supset S$ such that $D_{log}:=S_{log}-S$ is a simple normal crossing  divisor and the holomorphic symplectic form $\omega^{2,0}$ has poles of order 1 on the smooth part of $D_{log}$.\footnote{Partial log-compactifications are non-unique. A choice corresponds to the choice of irregular behavior at marked points in the case of connections and Higgs fields.}

With these data and a finite subset $R\subset D_{log}$ one can associate a family of $A_\infty$-categories\footnote{In fact it is an analytic family of $A_\infty$-categories, see~\cite{KontsevichSoibelmanTwistor} for details.}  $\cC_\zeta$, for $\zeta\in\P^1_\C$, in the following way. The fiber at $\zeta\ne0$ is the appropriately defined Fukaya category $\cF_\zeta(S,R)$ associated with the symplectic form $\omega_\zeta=\frac{\omega^{2,0}}{\zeta}+\zeta\overline{\omega^{2,0}}$, which naturally splits as a sum $Re(\omega_\zeta)+iB_\zeta$. We treat
\[
    B_\zeta=Im\left(\frac{\omega^{2,0}}{\zeta}+\zeta\overline{\omega^{2,0}}\right)
\]
as the $B$-field. The set $R$ is the finite subset of points where the closure of the support of an object intersects $D_{log}$.

The fiber at $\zeta=0$ is the ($A_\infty$-version of the) category $Coh^{ss,0}_c(S, R)$ of compactly supported  parabolic semistable coherent sheaves on $S$ of parabolic degree zero with pure 1-dimensional support satisfying certain conditions. In particular, the support is a curve whose closure intersects $D_{log}$ at the finite subset $R$. The details of the definitions of both categories, including the notion of parabolic degree and the precise conditions on the supports of objects near $D_{log}$, are worked out in~\cite{KontsevichSoibelmanTwistor}. Taking the inductive limit over all $R$ we obtain the categories $\cF_\zeta:=\cF_\zeta(S)$ and ${Coh}^{ss,0}_c(S)$, respectively.

The  category $Coh^{ss,0}_c(S)$ becomes the category of parabolic Higgs bundles of degree zero in the case $S=T^\ast X$, where $X$ is a curve, while $\cF_\zeta$ is equivalent in this case to the category of filtered local systems.
Finally, at  $\zeta=\infty$ we place the ``Hermitian dual'' category $(Coh^{ss,0}_c(S))^\ast=\overline{(Coh^{ss,0}_c(S))^{op}}$. The latter is defined as the opposite  category  with the opposite complex structure. This is our family $\cC_\zeta$.

\begin{definition}[see~\cite{KontsevichSoibelmanTwistor}] We say that the (analytic) family of $A_\infty$-categories $\cC_\zeta,\zeta\in\P_\C^1$ is a twistor family, if the following conditions are satisfied:

1) There is an equivalence of $A_\infty$-categories $G_\zeta: \cC_\zeta\simeq\cC^{\ast}_{-1/\overline\zeta}$,
where $\cC^{\ast}=\overline{\cC^{op}}$ is the Hermitian dual category.

2) $Hom$-spaces of each category are endowed with positive Hermitian forms compatible with the functors $G_\zeta$.
\end{definition}

It is explained in~\cite{KontsevichSoibelmanTwistor} that the above-defined family is a twistor family of categories. It gives rise to the analytic family of the corresponding stacks of objects. The following informal version of the \emph{generalized $P=W$ conjecture\/} was proposed by the third author (Y.S.) at the Workshop on Hall Algebras, Enumerative Invariants and Gauge Theories (Fields Institute, 2016).

\begin{conjecture}
The perverse filtration arising on the (properly understood) Borel-Moore homology of the moduli stack of objects $Coh^{ss,0}_c(S)$  coincides  with the weight filtration on the (properly understood) Borel-Moore homology of the moduli stack of  objects of  the Fukaya category $\cF_\zeta$ at $\zeta=1$.
\end{conjecture}

A more precise claim at the level of moduli spaces agrees with the conventional $P=W$ conjecture in the case $S=T^\ast X$.  One the other hand  the version of $P=W$ conjecture for $q$-difference equations (this corresponds to the complex surface $(\mathbb{C}^\ast)^2$) seems to be new.  The paper ~\cite{Davison2024HodgeStacks} is devoted to the discussion of the precise statement of the  version of the above \emph{stacky\/} conjecture in the case $S=T^\ast X$ (that is, for the moduli stack of Higgs bundles). One can hope that other results from loc. cit. admit upgrades in the framework of generalized non-abelian Hodge theory.

The question about motivic classes of the stacks of objects of the above-mentioned twistor family of categories, as well as the question about motivic generalized non-abelian Hodge theory, are widely open.

\subsection{Counting nilpotent pairs: double affine Grassmannians and geometric Satake correspondence}
Recall that the computation of the motivic classes of the stacks of nilpotent pairs is a key step in our computation (see~\cite[Sect.~4]{FedorovSoibelmansParabolic} and~\cite[Sect.~4]{FedorovSoibelmansIrregular}). In particular, we need to calculate the motivic class of the moduli stack $\Pair^{nilp,-}(X,D,\lambda)$ that parameterizes pairs $(\bE,\Psi\colon\bE\to\bE)$ , where $\bE$ is a level $D$ parabolic bundle and $\Psi$ is a nilpotent endomorphism of $\bE$ whose generic type is given by the partition $\lambda$. One of the main steps in our calculation is the following factorization formula for the graded motivic classes (see~\cite[Thm.~4.4]{FedorovSoibelmans})
\[
 [\Pair^{nilp,-}(X,D,\lambda)]=[\Pair^{nilp,-}(X,\emptyset,\lambda)]
 \prod_{x\in Supp(D)}V_x.
\]
The classes $V_x$ do not depend on the curve $X$, so they are local in nature. In fact, they can be described in terms of bundles on the formal disk. One cannot compute them naively using the techniques of Hall algebras, in particular, because the automorphism groups of points are infinite-dimensional.

We propose to replace Hall algebras by the spherical affine Hecke algebras of the relevant 2-dimensional versions of the affine Grassmannians. Then, symmetric functions will appear when we apply a version of the geometric Satake correspondence to the latter.

Let us explain it in the simpler case of nilpotent pairs over $\Spec\C$.That is, for nilpotent representations of the Jordan quiver (one vertex and one loop).

Let $Gr(n)=GL(n,\C((t)))/GL(n,\C[[t]])$ be the affine Grassmannian of $GL(n)$. There is a stratification $Gr(n)=\bigcup_\lambda Gr^\lambda(n)$, where $Gr^\lambda(n)$ is the $GL(n,\C[[t]])$-orbit of the diagonal matrix whose diagonal elements are $t^{\lambda_i}$, $\lambda=(\lambda_1\ge\lambda_2\ge\ldots)$. Let $Gr^+(n):=\bigcup_\lambda Gr^\lambda(n)$ be the positive Grassmannian. Note that $Gr^+(n)$ classifies lattices $L\subset\C[[t]]^n$, that is, $\C[[t]]$-submodules of finite co-rank over $\C$. Then, $t$ induces a nilpotent endomorphism of $\C[[t]]^n/L$. Denote the quotients stacks $GL(n,\C)\backslash Gr^+(n)\subset GL(n,\C[[t]])\backslash Gr(n)$ by $Mod^+(n)\subset Mod(n)$.

Let $R_d$ denote the moduli stack $GL(d)\backslash Nil_d$, where $Nil_d$ is the nilpotent cone in the Lie algebra $\gl(d,\C)$. In other words, this is the moduli stack of $d$-dimensional representations of the Jordan quiver. Set $R:=\bigsqcup_{d\ge0} R_d$. We obtain a morphism of stacks
\[
    F\colon Mod^+(n)\to R\colon (L\subset\C[[t]]^n)\mapsto(\C[[t]]^n/L,\bar t),
\]
where $\bar t$ is the automorphism of $\C[[t]]^n/L$ induced by $t$. One shows that the pullback along this morphism induces a surjective homomorphism of the Hall algebra of the category of nilpotent representations of the Jordan quiver to the spherical affine Hecke algebra of $GL(n)$. This surjective homomorphisms, in turn, induces an isomorphism of the Hall algebra and the projective limit as $n\to \infty$ of the affine spherical Hecke algebras.\footnote{The latter can be thought of as the spherical affine Hecke algebra of the inductive limit $Mod^+:=\varinjlim_n Mod^+(n)=\bigcup_n Mod^+(n)$, where the morphism $Mod^+(n)\to Mod^+(m)$ is induced be an embedding $\C[[t]]^n\hookrightarrow\C[[t]]^m$ when $n\le m$.}

The nilpotent orbit corresponding to the partition $\lambda$ corresponds to $Gr^\lambda(n)$. One can show that the delta function of the nilpotent orbit goes to the Hall--Littlewood polynomial after the Satake isomorphism, which identifies the spherical affine Hecke algebra with the algebra of symmetric functions. It might be interesting to compare this approach with the one from~\cite{Chaudouard2019Counting}.

Now, the idea for an alternative approach to the computation of the local factors $V_x$ is to replace nilpotent representations of the Jordan quiver in vector spaces by those in free $\C[[t]]$-modules, and to replace the positive affine Grassmannian by its appropriate 2-dimensional version.


\newcommand{\etalchar}[1]{$^{#1}$}

\end{document}